\documentclass{amsart}

\newtheorem{Def}{Definition}[section]

\newtheorem{Cor}[Def]{Corollary}
\newtheorem{Thm}[Def]{Theorem}
\newtheorem{Lem}[Def]{Lemma}
\newtheorem{Rem}[Def]{Remark}
\newcommand{\Cplx}{\mathbb{C}}
\newcommand{\Real}{\mathbb{R}}
\newcommand{\Zahl}{\mathbb{Z}}

\newcommand{\Nash}{{X'}}

\begin{document}

\title{Pseudo vector bundles and quasifibrations}
\author{Martin A. Guest, Micha\l \ Kwieci\'nski and Boon-Wee Ong }
\address{
{\it M.A.G., M.K.:}  Department of Mathematics,
Tokyo Metropolitan University,
Minami - Ohsawa 1-1, Hachioji-shi, Tokyo 192-0397, Japan
}
\address{
({\it M.A.G. is on leave from:} Department of Mathematics, University of
Rochester, Rochester, NY 14627, USA;
{\it M.K. is on leave from:}
Uniwersytet Jagiello\'nski, Instytut Matematyki,
\hbox{ ul. Reymonta 4}, 30-059 Krak\'ow, Poland.)
}
\address{
{\it B.O.:} Department of Mathematics, University of
Rochester, Rochester, NY 14627, USA
}
\address{\phantom{a}}
\email{M.A.G.: martin@comp.metro-u.ac.jp; M.K.:
michal@math.metro-u.ac.jp, kwiecins@math.metro-u.ac.jp;
B.O.: amos@math.rochester.edu}
\thanks{M.K. is supported by the Japan Society for the Promotion of Science
and a Monbusho grant.}


\begin{abstract}
We prove a topological result concerning the kernel $\ker d$ of
a morphism $d:E\to F$ of holomorphic vector bundles over a
complex analytic space.  As a consequence, we show that the
projectivization $\Bbb P(\ker d)$ is a quasifibration up
to some dimension. We give
an application to the Abel - Jacobi map of a Riemann surface,
and to the space of rational curves in the symmetric product
of a Riemann surface.

\end{abstract}

\maketitle

\section{Introduction}

Let $d:E\to F$ be a morphism of (smooth) vector bundles over
a manifold $M$.  If the rank of $d_m:E_m\to F_m$ is
independent of $m$, the space
$$
L=\ker d = \bigcup_{m\in M} \ker  d_m
$$
is a subbundle of $E$.  In general, however, the fibres of
the map $\pi:L \to M$ are vector spaces of
varying dimensions, and the topological behaviour of
$\pi$ can be very complicated; it will certainly not be locally trivial.

In the complex analytic category, the topological behaviour
of $\pi$ is more predictable.  We shall show that $\pi$
\lq\lq resembles" a vector bundle  from the viewpoint of homotopy theory,
and in particular that the projectivization
$\Bbb P(\pi): \Bbb P (L) \to M$ is a quasifibration
up to some dimension (Theorem \ref{projectiv}).

Our results are valid in the more general situation
where $d:E\to F$ is a morphism of holomorphic vector
bundles over a complex analytic space $X$, and we shall always work in
this generality from now on, also assuming that $X$ is nonempty.
The space $L$
is an analytic subspace of the total space of
$E$, and we shall call the map
$\pi:L \to X$ a {\it pseudo vector bundle}.
A classical example of a pseudo vector bundle is the Zariski tangent space
(to a singular complex algebraic variety).
Pseudo vector bundles have been studied under the label
\lq\lq linear spaces" in analytic geometry
(see \cite{Fi}, \cite{R}, \cite{S}) and as (spectra of)
\lq\lq symmetric algebras" in
commutative algebra (see \cite{HR}, \cite{V}).
There are still many open problems concerning these
objects ---
for example finding necessary and sufficient conditions for their
irreducibility or equidimensionality (see \cite{Kw}, \cite{V}).
It should be stressed that $L$
is conceptually different from the kernel sheaf of the induced map
of sheaves of sections of $E$ and $F$. Indeed, the induced map may
be injective in cases when the map of vector bundles is not.

Theorem \ref{filtr} below is our basic result on the topological
structure of pseudo vector bundles.
Before stating it, let us just remark
that the concepts of pullback, restriction and subbundle apply
to pseudo vector bundles in the obvious way.

\begin{Thm}
\label{filtr}
There exists a (unique) filtration of $X$ by closed analytic subsets :
$$X=X_k\supset X_{k+1}\supset\cdots\supset X_l=\emptyset\ ,$$
such that $X_k\neq X_{k+1}$ and
\begin{itemize}
\item
for each $i$, either $X_i\setminus X_{i+1}$ is empty or the restriction
$L\vert_{X_i\setminus X_{i+1}}$
is a vector bundle of rank $i$.
\end{itemize}
Furthermore, for each $i$ there exists
\begin{itemize}
\item an open neighbourhood $U_i$ of $X_{i+1}$ in $X_i$,
\item a (strong) deformation retraction $r_i:U_i\to X_{i+1}$,
\item a (strong) deformation retraction
$\rho_i:L\vert_{U_i}\to L\vert_{X_{i+1}}$,
\end{itemize}
such that
\begin{itemize}
\item $r_i^t\circ\pi=
\pi\circ\rho_i^t$ on $L\vert_{U_{i}}$,
where $r_i^t,\rho_i^t$ (respectively) are homotopies of $r_i,\rho_i$
to the identity maps of $U_i,L\vert_{U_i}$.
\item $\rho_i$ restricted to each fibre of $L$ is a linear
monomorphism.
\end{itemize}
\end{Thm}

The first statement of the
theorem is easily proved by considering a local model
of $d$, that is a matrix $A$, whose entries are holomorphic functions
on an open subset of $X$. Then the sets $X_i$ of the
filtration are described as the zero sets of ideals generated by
minors (of
the appropriate size) of $A$ and hence are automatically analytic.
These ideals are known as Fitting ideals and have been a subject
of intensive study in commutative algebra (see e.g. \cite{Eisenbook}).
The remaining statements will be proved later in the paper.

\begin{Def}
\label{rank}
We call the integer $k$ appearing in  Theorem \ref{filtr} the 
{\it minimal rank} of $L$.
If $X_k\setminus X_{k+1}$ is dense in $X$, then we
call $k$ the {\it generic rank} of $L$ and say that $L$ has a generic rank.
\end{Def}
\noindent Of course, if $X$ is irreducible then $L$ always has a generic rank.

It should be noted that a pseudo vector
bundle is not necessarily a \lq\lq direct sum" of
vector bundles on analytic subsets of $X$. Indeed, even
in the complex analytic situation, the behaviour of
the fibres of $L\vert_{X_i\setminus X_{i+1}}$ near $X_{i+1}$ may
be rather complicated. It is in general impossible to extend
the vector bundle $L\vert_{X_i\setminus X_{i+1}}$ to a vector
bundle $K$ on $X_i$ such that $K$
is a subbundle of $L$. The simplest example
of such a situation is the kernel of the morphism
from a rank two trivial bundle
to a rank one trivial bundle
over $\Cplx^2$ given by the matrix $(x,y)$.

It is an immediate consequence of Theorem \ref{filtr} and the
Dold-Thom criterion \cite{Dold-Thom} that $\Bbb P(\pi)$ is a
quasifibration up to dimension $2k-1$;
this statement appears later as Theorem \ref{projectiv}.
Quasifibrations have played an important but narrowly focused
role in algebraic topology (see for example \cite{Dold-Thom},
\cite{Mc}, \cite{Mc-SLN}, \cite{Mc-Se}, \cite{Se}, \cite{Ag-Pr}).
Our result gives a simple but very general family of new examples.

We shall give a brief review of the theory of quasifibrations
in section 2, and then in section 3 we give the proof
of Theorem \ref{filtr}. In section 4, we give
an application to the Abel - Jacobi map of a Riemann surface,
and to spaces of rational curves in symmetric products
of Riemann surfaces.
Finally, in section 5, we mention an example related to 
the theory of harmonic maps, and also a
counterexample to the statement of Theorem \ref{filtr}
in the category of smooth vector bundles.

\section{Quasifibrations}

We begin by recalling the standard fact from homotopy theory
that \lq\lq any continuous map is homotopic to a fibration".
Let $f:X\to Y$ be a continuous map of topological spaces.
The {\it homotopy fibre} of $f$ over $y$ is defined to be
the space
$$
H_y = \{ (x,\gamma) \ \vert\
x\in X, \gamma:[0,1]\to Y, \gamma(0)=y, \gamma(1)=f(x) \},
$$
i.e. (continuous) paths in $Y$ from $y$ to $f(x)$.
This is the fibre over $y$ of the fibration
$$
\tilde f:
\{ (x,\gamma) \ \vert\
x\in X, \gamma:[0,1]\to Y, \gamma(1)=f(x) \} \to Y,
\quad (x,\gamma) \mapsto \gamma(0).
$$
Since the domain of $\tilde f$ is homotopy equivalent
to $X$, it follows that --- \lq\lq up to homotopy" --- we
may identify $f$ with the fibration $\tilde f$.

There is a natural inclusion map $f^{-1}(y)\to H_y$.  If this is a homotopy
equivalence for all $y\in Y$, we say that $f$ is a
{\it quasifibration.}  For reasonable spaces $X$ and
$Y$, any fibration has this property, i.e. any fibration
is a quasifibration.  But in general the property of being a
quasifibration is weaker than the property of being a
fibration.

It is well known that there is a long exact sequence of
homotopy groups
$$
\dots \to \pi_i f^{-1}(y) \to \pi_i X \to \pi_i Y \to
\pi_{i-1} f^{-1}(y) \to \dots
$$
if $f$ is a fibration.  This remains true if $f$ is a
quasifibration, since we may replace $H_y$ by $f^{-1}(y)$ in the long
exact sequence of the fibration $\tilde f$.

In this paper we shall need the following weaker concept:

\begin{Def}
\label{quasifib-k}
A map $f:X\to Y$ is a quasifibration up to dimension
$p$ if the inclusion $f^{-1}(y)\to H_y$ is a (weak) homotopy
equivalence up to dimension $p$, i.e. if this inclusion map
induces isomorphisms of homotopy groups $\pi_i$ for
$i<p$, and a surjection for $i=p$.
\end{Def}

Quasifibrations first appeared in \cite{Dold-Thom}.
More generally, the concept of
\lq\lq homology fibration (up to dimension $p$)"
was introduced in \cite{Mc} (Lemma 4.1) and \cite{Mc-Se}.
In all of these papers the main purpose was to study
maps between configuration spaces.  However, there
are other non-trivial applications, for example an
approach via quasifibrations to the Bott Periodicity
Theorem was suggested in \cite{Mc-SLN} and carried out in 
\cite{Ag-Pr}, and homology fibrations
(up to dimension $p$) were used to study spaces of
rational functions in \cite{Se}.

To prove that a map is a quasifibration,
the \lq\lq Dold-Thom Criterion" of \cite{Dold-Thom}
is often used.  This criterion may be modified
in an obvious way to prove that a map is a quasifibration
up to dimension $p$. The following weaker version of the
result of \cite{Dold-Thom} will be sufficient for our purposes.

\begin{Thm}
\label{DTC}
Let $f:Y\to X$ be a map, $p$ a positive integer.
Assume that $X$ has a filtration by closed subsets :
$$X=X_k\supset X_{k+1}\supset\cdots\supset X_l=\emptyset\ ,$$
such that
\begin{itemize}
\item
for each $i$,  the restriction of $f$ to
$f^{-1}({X_i\setminus X_{i+1}})$
is a fibration.
\end{itemize}
Assume further that for each $i$ there exists
\begin{itemize}
\item an open neighbourhood $U_i$ of $X_{i+1}$ in $X_i$,
\item a (strong) deformation retraction
$r_i:U_i\to X_{i+1}$,
\item a (strong) deformation retraction
$\rho_i:f^{-1}(U_i)\to f^{-1}(X_{i+1})$,
\end{itemize}
such that
\begin{itemize}
\item $r_i^t\circ f=
f\circ\rho_i^t$ on $f^{-1}(U_i)$,
where $r_i^t,\rho_i^t$ (respectively) are homotopies of $r_i,\rho_i$
to the identity maps of $U_i,f^{-1}(U_i)$.
\item $\rho_i$ restricted to each fibre of $f$ is a homotopy
equivalence up to dimension $p$.
\end{itemize}
Then $f$ is a quasifibration up to dimension $p$.
\end{Thm}

Combining this with Theorem \ref{filtr}, we obtain the statement in the
introduction concerning a projectivized pseudo vector bundle
$\Bbb P(\pi)$:

\begin{Thm}
\label{projectiv}
Let $k=\text{min}\{ \text{dim}_{\Cplx}
\ker  d_x \ \vert\ x\in X \}$.
Then the map
$$
\Bbb P(\pi): \Bbb P(L) = \Bbb P(\ker  d) \to X
$$
is a quasifibration up to dimension $2k-1$.
\end{Thm}

\begin{proof}
By Theorem \ref{filtr}, the \lq\lq attaching map"
$\Bbb P(\rho_i \vert_{\ker  d_x})$ is a linear
inclusion of the form $\Cplx P^s \to \Cplx P^t$,
with $k-1 \le s \le t$. It is
well known that any such map is a homotopy equivalence
up to dimension $2s+1$ (for example, because $\Cplx P^t$
may be constructed topologically from $\Cplx P^s$ by
adjoining cells of dimensions $2s+2$, $2s+4$, \dots, $2t$).
\end{proof}

\section{Pseudo vector bundles.}

This section is devoted to the proof of the remaining parts
of Theorem \ref{filtr}.
For simplicity, we shall write $V=X_k\setminus X_{k+1}$ and
$S=X_{k+1}$ ($S$ for \lq\lq singular set").

It is well known that one can simplify the structure of $L$ by
an appropriate blowing-up procedure (see \cite{R} or \cite{S}), known as the
Nash modification. For our purposes we shall need the following result.

\begin{Lem}
\label{Nash}
Suppose $L$ has a generic rank. Then there exists an analytic space
$\Nash$ and a proper surjective analytic map $\nu:\Nash \to X$,
such that
\begin{itemize}
\item the inverse image $\nu^{-1}V$ is dense in $\Nash$,
\item the restriction $\nu\vert_{\nu^{-1}V}:\nu^{-1}V\to V$ is
a biholomorphism,
\item there exists a vector bundle $L'$, which is a subbundle of
the pseudo bundle $\nu^*L$, such that
$$L'\vert_{\nu^{-1}V}\ =\ \nu^*L\vert_{\nu^{-1}V}\ .$$
\end{itemize}
\end{Lem}

Actually, $\nu$ becomes unique if we require an additional
minimality condition, which we do not need here.
The pair $(\Nash,\nu)$ may be constructed as follows.
The morphism $d$ defines a section of the Grassmannian bundle:
$$\sigma: V\to {\text{Grass}}_k(E)\ ,$$
$$x\to \ker d_x\ .$$
This section is defined over $V$ only. We define
$\Nash$ to be the closure of the image of $\sigma$ and
the mapping $\nu:\Nash\to X$ to be the restriction of the
natural projection ${\text {Grass}}_k(E)\to X$. All conditions
of Lemma \ref{Nash} are
now satisfied. In particular, the bundle $L'$ may be taken to be the
restriction of the tautological bundle on the Grassmannian bundle.
For the details, we refer to \cite{R}. We just mention that
the analyticity of $\Nash$ is proved by doing a simple calculation
in Pl\"ucker coordinates, which shows that $\Nash$ is the closure
of an analytically constructible set. \qed

To find the $i$-th neighbourhood and retractions in Theorem \ref{filtr},
one only cares about what goes on inside $X_i$. Therefore, one may assume
$i=k$ and the theorem will follow from the next statement.

\begin{Lem}
\label{one}
Let $\pi:L\to X$ be a pseudo vector bundle on an analytic space $X$
and let $S$ be an analytic subset of $X$, $V=X\setminus S$.
Assume, that $L\vert_{V}$ is a vector bundle of rank $k$.
Then there exists
\begin{itemize}
\item a neighbourhood $U$ of $S$ in $X$,
\item a (strong) deformation retraction $r:U\to S$,
\item a (strong) deformation retraction
$\rho:L\vert_{U}\to L\vert_{S}$,
\end{itemize}
such that
\begin{itemize}
\item $r^t\circ\pi=
\pi\circ\rho^t$ on $L\vert_{U}$,
where $r^t,\rho^t$ (respectively) are homotopies of $r,\rho$
to the identity maps of $U,L\vert_{U}$.
\item $\rho$ restricted to each fibre of $L$ is a linear
monomorphism.
\end{itemize}
\end{Lem}

\bigskip

\begin{proof}
We shall prove the above lemma in three steps.

\medskip

\noindent {\bf Step 1.} Lemma \ref{one} is true in the case when $L$ itself
is a vector bundle.

The existence of $U$ and $r$ is obvious (for example by
choosing compatible triangularizations of $X$ and $S$).
It is a standard fact concerning vector bundles
(see \cite{A}, Lemma 1.4.3) that we have the following
isomorphism (of $\mathcal C^0$ vector bundles)
$$L\vert_{U}\cong r^*(L\vert_{S}) .$$
Then $\rho$ can be constructed as the composition of
this isomorphism with the natural map
$r^*(L\vert_{S})\to L\vert_{S}$.

\medskip

\noindent {\bf Step 2.}  Lemma \ref{one} is true in the case when $L$ has a
generic rank.

The assumption allows us to use Lemma \ref{Nash} to produce
$\Nash$, $\nu$ and $L'$. Now, $L'$ is a true vector
bundle. Hence we can apply the result of Step 1
to the space $\Nash$ with the subset $S'=\nu^{-1}(S)$
and the vector bundle $L'$ to obtain $U'$, $r'$ and $\rho'$.
Next, consider the subbundle inclusion
$$\iota:\ L'\vert_{S'}\subset \nu^*L\vert_{S'}$$
and
$$\rho':L'\vert_{U'}\to L'\vert_{S'}\ .$$
Glue $\iota\circ\rho'$ (defined on $L'\vert_{U'}$)
with the identity on ${\nu^*L\vert_{S'}}$
to obtain a retraction
$$\rho'':
\nu^*L\vert_{U'}\to \nu^*L\vert_{S'}\ .$$
The map thus defined is continuous, because we have constructed it as
the glueing of maps  on two closed sets, which agree
on their intersection ($L'\vert_{S'}$).

We now have suitable retractions on the pullback of $L$ by $\nu$.
The remaining problem is to push them forward back again.

First of all, since $\nu $ is proper and surjective,
by elementary analytic topology one proves that
the image of $U'$ contains a neighbourhood of
$S$ (we leave this to the reader).
Since $\nu$ is a homeomorphism on the complement of $S'$,
the image
$U=\nu(U')$ is actually an open neighbourhood of $S$.
Shrinking $U$ if necessary, we may suppose that $U'=\nu^{-1}(U)$.
Then $\nu\vert_{U'}:U'\to U$ is also proper.

It is now easy to find the retraction $r$. Notice that
$\nu\vert_{S'}\circ r'$ is constant on each fibre of $\nu\vert_{U'} $.
Since $\nu\vert_{U'}:U'\to U$ is surjective, this implies the existence of
a unique map $r:U\to S$, such that
$$r\circ\nu\vert_{U'}=\nu\vert_{S'}\circ r'\ .$$
The continuity of $r$ easily follows from the fact that
$\nu \vert_{U'}:U'\to U$ is proper and surjective and the other properties
are immediately transferred from those of $r'$.

The construction of $\rho$ from $\rho''$ follows along the same lines.
Instead of $\nu $, we use the induced map
$\tilde\nu:\nu^*L\to L$ and remark that it is proper and
surjective as the base change (fibred product)
of the proper surjective map $\nu$
by the projection $\pi:L\to X$. Again, $\rho$ is the unique map
which satisfies
$\rho\circ\tilde\nu\vert_{\nu^*L\vert_{U'}}=
{\tilde\nu}\vert_{\nu^*L\vert_{S'}}\circ\rho''$.

\medskip

\noindent {\bf Step 3.} The general case.

Let $\tilde S=\bar V\cap S$.
Now $\bar V$ (the closure of $V$) is an analytic set, as the closure
of a constructible analytic set.
We can now apply the result of the previous step
to $\bar V\supset \tilde S$ and the restriction of $L$.
Taking the union of the neighbourhood with $S$ and
extending the retractions
by the identity on $S$ and on $L\vert_S$ one obtains
the desired neighbourhood and retractions in the general case.
\end{proof}

With a little more effort (essentially Lemma \ref{one} applied inductively)
one can prove a stronger result, which we do not need here, but which
is perhaps worth noticing :

\begin{Rem}
Lemma \ref{one} is true without the assumption that $L\vert_{V}$ is a
vector bundle.
\end{Rem}

\section{The Abel - Jacobi Map.}

Let $M$ be a compact Riemann surface of genus $g$, and let
$J(M)$ be the Jacobian variety of $M$ (a complex torus
of dimension $g$). Let $\text{Sp}^d(M)$ be the
$d$-th symmetric product of $M$, i.e.
$\text{Sp}^d(M) = M^d / \Sigma_d$ where $\Sigma_d$ is
the symmetric group.  Although the action of $\Sigma_d$
is not free, it is known that $\text{Sp}^d(M)$ has the
structure of a complex manifold of dimension $d$
(see \cite{Gunning}, \S 3).

The Abel-Jacobi map $j:\text{Sp}^d(M) \to J(M)$ is
a holomorphic map, which may be described as follows:
if $\{ m_1,\dots,m_d\}$ is an element of
$\text{Sp}^d(M)$, then $j(\{ m_1,\dots,m_d\})$
is the isomorphism class of the holomorphic line bundle
$L\otimes L_0$ on $M$, where $L$ is a line bundle
corresponding to the divisor $\sum_{i=1}^d m_i$ and $L_0$
is a fixed line bundle of degree $-d$.

It is a result of Mattuck that $j$ is a holomorphic
fibre bundle, with fibre $\Cplx P^{d-g}$, if $d\ge 2g-1$.
A proof of Mattuck's theorem, together with a more
explicit description of $j$, may be found in \cite{Gunning}
(see also \cite{Ong}).
It can be shown from this description that $j$ is in fact a
projectivized pseudo vector bundle, for any $d\ge g$. (If $d<g$,
$j$ cannot be surjective.) Moreover, for $g\le d\le 2g-2$,
this pseudo vector bundle has generic rank $d-g+1$.
Hence we obtain:

\begin{Thm}
\label{Amos}
For $g\le d\le 2g-2$, the Abel-Jacobi map
$j:\text{Sp}^d(M) \to J(M)$
is a quasifibration up to dimension $2(d-g)+1$.
\end{Thm}

\noindent This result may also be verified by a direct
(but unilluminating) calculation
since the homotopy types of $\text{Sp}^d(M)$, $J(M)$,
and the fibres of $j$ are all well-known (see \cite{Ong}).

Let $\text{Hol}_k(S^2,\text{Sp}^d(M))$ denote the space
of holomorphic maps $f:S^2\to \text{Sp}^d(M)$ whose homotopy class
$[f]\in \pi_2 \text{Sp}^d(M) \cong \Zahl$ is $k$ ($\ge 0$).
Taking the induced topology from the corresponding space
$\text{Map}_k(S^2,\text{Sp}^d(M))$ of continuous maps, we
have a continuous inclusion map $i_k:\text{Hol}_k(S^2,\text{Sp}^d(M))
\to \text{Map}_k(S^2,\text{Sp}^d(M))$.  There are many examples
of compact complex manifolds $X$ for which the analogous inclusion
map $\text{Hol}(S^2,X)\to\text{Map}(S^2,X)$ is a homotopy equivalence
up to some dimension (where the dimension depends on the homotopy class).
A detailed explanation of this type of result and its significance can be found,
for example, in \cite{Gu} and \cite{Hu}.  We shall prove
that the space $X=\text{Sp}^d(M)$ is another such example, 
in fact one which is rather different from those considered
so far in the literature.  The main tool is the following
extension of Theorem \ref{Amos}:

\begin{Thm}
\label{Hol}
For $g\le d\le 2g-2$, the  map
$j_k:\text{Hol}_k(S^2,\text{Sp}^d(M))\to J(M)$,
defined by $j_k(f)=j\circ f$,
is a quasifibration up to dimension $2(d-g)-1$.
\end{Thm}

\begin{proof}
The map is well defined as any holomorphic map $j\circ f:S^2\to J(M)$ is necessarily
constant.  Hence, the image of $f$ lies entirely in the fibre of  $j$ over
the point $j\circ f$. To prove the theorem,
exactly the same argument may be used as in the case $k=0$,
except that for $k\ge 1$ we need to know that, if $d-g\le s\le t$, the inclusion
$\text{Hol}_k(S^2,\Cplx P^s) \to \text{Hol}_k(S^2,\Cplx P^t)$ is a homotopy
equivalence up to dimension at least $2(d-g)-1$.  This follows (for example)
from Segal's theorem (\cite{Se}) that the inclusion
$\text{Hol}_k(S^2,\Cplx P^N) \to \text{Map}_k(S^2,\Cplx P^N)$ is
a homotopy equivalence up to dimension $(2N-1)k$, because the
map $\text{Map}_k(S^2,\Cplx P^s) \to \text{Map}_k(S^2,\Cplx P^t)$ is
a homotopy equivalence up to dimension at least $2(d-g)-1$ (by the argument
in the proof of Theorem \ref{projectiv}).
\end{proof}

\begin{Cor}
\label{Amos2}
For $g\le d\le 2g-2$ and $k\ge 1$,
the natural inclusion map
\newline
$i_k:\text{Hol}_k(S^2,\text{Sp}^d(M))
\to \text{Map}_k(S^2,\text{Sp}^d(M))$
is a homotopy equivalence up to dimension $2(d-g)-1$.
\end{Cor}

\begin{proof} From Theorem \ref{Amos}, 
the map $j^\prime_k:\text{Map}_k(S^2,\text{Sp}^d(M))\to 
\text{Map}_k(S^2,J(M))$,
defined by $j^\prime_k(f)=j\circ f$,
is a quasifibration  up to dimension $2(d-g)-1$.  The space
$\text{Map}_k(S^2,J(M))$ is homotopy equivalent to $J(M)$ itself,
so we may compare the long exact sequences of homotopy groups
of $j_k$ and $j^\prime_k$; the desired result then follows
from the Five Lemma and the result of
Segal quoted in the proof of Theorem \ref{Hol}.
\end{proof}

For $d>2g-2$ all relevant maps in the proof of the corollary
are quasifibrations, and
so $i_k$ is a homotopy equivalence up to dimension $[2(d-g)-1]k$ 
in this case.  Thus, the \lq\lq approximation\rq\rq\  improves as $k\to\infty$,
as in all other previously known examples of this type.  But in the
range $g\le d\le 2g-2$ no such improvement appears to be possible.

Corollary \ref{Amos2} may be relevant to the quantum cohomology
of $\text{Sp}^d(M)$ --- see the recent article \cite{Be-Th}.

\section{Further examples}

Many examples of pseudo vector
bundles may be constructed explicitly by writing down matrix-valued
complex analytic functions. It is less
easy to find examples where the total space is a smooth
variety --- the Abel-Jacobi map being one such case.  Another
case arises naturally in connection with the theory of
harmonic maps in differential geometry.  It is well known that
all harmonic maps from $\Cplx P^1$ to  $\Cplx P^2$ may be
constructed from holomorphic maps , and it was proved in
\cite{Cr} that each component of the space
of harmonic maps  $\Cplx P^1 \to \Cplx P^2$ is a smooth variety. 
The description of this space makes use of a pseudo vector
bundle whose total space is smooth; this is described in detail
in section 3 of \cite{Le-Wo}.

Finally, we note that Theorem \ref{filtr} is {\it false}
in the smooth category.  There is no guarantee that
a suitable neighbourhood $U_i$ of $X_i$ can found; for
example the $1\times1$ matrix-valued function
 \[
        d(u) = \left\{
        \begin{array}{ll}
                e^{- \frac{1}{u^2}}   \sin \frac1u
                & \text{if $u \neq 0$} \\
                0     &  \text{if $u = 0$}
        \end{array}    \right.
 \]
has rank zero on the subset
\[
        X = \left\{  u = 0 \ \text{or}\  \frac{1}{\pi n},
                        \quad n \in \Zahl\setminus\{0\} \right\}
\]
of $\Real$, and this set admits no neighbourhood
of which it is a deformation retract.  Even when
a suitable neighbourhood $U_i$ of $X_i$ exists, the
following example shows that the fibres of the pseudo
vector bundle can behave badly near $X$.
Consider the smooth map
$$
d:M\times \Real^2 \to M\times \Real,
\quad
d(m,(x,y)) = (m,x f(m) \cos(m^{-1}) + y f(m) \sin(m^{-1}) )
$$
where $M=\Real$, and where $f$ is a smooth function which
vanishes to all orders at zero. Over $\{ 0\}$ the fibre of $\ker d$ is
$\Real^2$, and over $M \setminus \{ 0\}$ the space $\ker d$
is a locally trivial bundle with fibre $\Real$. But clearly
the conclusions of Theorem \ref{filtr} do not hold.

\end{document}